\documentclass{article}
\usepackage{amssymb}
\usepackage{amsfonts}
\usepackage{amsmath}
\usepackage{color}
\usepackage{graphicx}

\newtheorem{lemma}{Lemma}
\newtheorem{theorem}{Theorem}

\newtheorem{problem}{Problem}
\newtheorem{conjecture}{Conjecture}
\newcommand{\DONOTTEX}[1]{}

\newcommand{\NEWTEXT}[1]{{#1}} \newcommand{\OLDTEXT}[1]{}

\DeclareMathOperator{\ex}{ex}

\begin{document}
\title{{\bf Maximal ambiguously \\ \boldmath$k$-colorable graphs}}
\author{Matthias Kriesell}

\maketitle

\begin{abstract}
  \setlength{\parindent}{0em}
  \setlength{\parskip}{1.5ex}

  A graph is {\em ambiguously $k$-colorable} if its vertex set admits two distinct partitions each into at most $k$ anticliques.
  We give a full characterization of the maximal ambiguously $k$-colorable graphs in terms of $k \times k$-matrices.
  As an application, we calculate the maximum number of edges an ambiguously $k$-colorable graph can have,
  and characterize the extremal graphs.
        

  
\end{abstract}

\maketitle

\section{Introduction}

An {\em anticlique} of a graph $G$ is a set of pairwise nonadjacent vertices of $G$, and a {\em $k$-coloring} of $G$ is
a partition of $V(G)$ into at most $k$ anticliques. Graphs with at least one $k$-coloring are {\em $k$-colorable}, and we call those
with more than one $k$-coloring {\em ambiguously $k$-colorable}. 
A graph is {\em maximal ambiguously $k$-colorable}
if it is ambiguously $k$-colorable but adding any edge between distinct nonadjacent vertices produces a graph which is not.
We give a full description of the maximal ambiguously $k$-colorable graphs in terms of quadratic matrices.

\NEWTEXT{
The initial motivation of studying these graphs have been some results and observations in \cite{Kriesell2016} on the conjectures of {\sc Hadwiger} and {\sc Seymour}
restricted to the case of {\em uniquely $k$-colorable graphs}, i.~e. graphs admitting only one coloring. {\sc Hadwiger}'s Conjecture \cite{Hadwiger1943} states
that, for every $k$, every graph either has a $(k-1)$-coloring or admits a set of $k$ many nonempty, connected, pairwise disjoint
and pairwise adjacent subgraphs (a so-called {\em clique minor of order $k$}). The conjecture has been verified for $k \leq 6$,
where for each of $k=5$ and $k=6$ it turned out to be equivalent to the $4$-color-theorem (see \cite{RobertsonSeymourThomas1993}). 
Since a uniquely $k$-colorable graph is not $(k-1)$-colorable
unless it is a complete graph on less than $k$ vertices, it follows, for $k \leq 6$, that a uniquely $k$-colorable graph $G$ is a complete graph
on less than $k$ vertices or has a clique minor of order $k$; in \cite{Kriesell2016} this has been proven {\em without using the $4$-color-theorem}
(or any other result for which we know computer-aided proofs only). For {\sc Seymour}'s Conjecture (see \cite{Blasiak2007}) that every graph on $n$ vertices without antitriangles
(i.~e.~anticliques of order $3$) admits a set of at least $n/2$ many pairwise disjoint pairwise adjacent complete subgraphs of order $1$ or $2$ (a so-called {\em shallow clique minor}),
the situation is even better:
A uniquely $k$-colorable graph $G$ without antitriangles is either a complete graph on less than $k$ vertices or it has a shallow clique minor of order $k$
(where $k \leq |V(G)|/2$ as $G$ contains no antitriangles).

These two results suggest that {\sc Hadwiger}'s conjecture and its relatives are substantially easier to deal with if
the number of $k$-colorings is just one (or, more general: limited).
What happens in the other case, where there are at least two (or, in general: at least $d$) $k$-colorings?
In order to understand this it is reasonable to look at the {\em extremal} case first, where we have as many edges as possible,
because this simplifies finding a sufficiently large clique minor. More generally, one could try to study the
{\em saturated} case, where the addition of an edge not present yet kills the property of having at least two (or, in general: $d$) $k$-colorings,
because it appears that in order to describe the extremal graphs one has to do this anyway.

For the case of maximal ambiguously $k$-colorable graphs both the conjectures of {\sc Hadwiger} and {\sc Seymour} are finally verified
by our main result. However, the argument involves only the (easy) observation that every such graph is an induced subgraph of
$\overline{K_k \times K_k}[K_\ell]$ for sufficently large $\ell$, as the latter graph is known to be perfect;
but no such argument will work for what is naturally be defined to be a maximal $d$-fold $k$-colorable graph, $d>2$ (see section \ref{SOpenquestions} for the details),
and so to understand the phenomena by which these graphs are ruled for $d=2$, as supported by the present paper,
may be very useful for the more general case.

There are several algebraic characterizations of unique $k$-colorability (see \cite{HillarWindfeldt2008} and the papers mentioned there in the introduction), and
it would be definitely very interesting (but is nonetheless postponed to the future) to relate them to our main result. However, as we are
lacking an ``algebraic counterpart on the minor side'', for example a characterization of graphs admitting a (perhaps very special) clique minor of order $k$,
it seems to be unlikely that there is an algebraic road leading to {\sc Hadwiger}'s Conjecture.
}

Let us describe our main theorem.
Let $A$ be a $k \times k$-matrix where all entries are non-negative integers.
$A$ is {\em tiny} if it is a diagonal matrix with exactly one entry $2$, all others at most $1$, and at least two diagonal entries $0$.
$A$ is {\em small} if it is a diagonal matrix with at least one entry $2$, all others at most $2$, and exactly one diagonal entry $0$.
$A$ is {\em special} if all diagonal entries are nonzero, exactly one off-diagonal entry is $1$, and all others are $0$.
\OLDTEXT{
$A$ is {\em normal} if it is a block diagonal matrix with diagonal blocks $M,D$,
where $D$ is a diagonal matrix with all diagonal entries nonzero, and $M$ has each of the following properties:
(i) All diagonal entries are nonzero,
(ii) $M$ is of size $r \geq 2$ and {\em fully indecomposable},
that is, it does not admit an $s \times (r-s)$ zero submatrix, where $s \in \{1,\dots,r-1\}$, and
(iii) whenever $M(i,j) \geq 2$ for some $i \not=j$, then there exists a sequence $f_0,\dots,f_\ell$ from $\{1,\dots,r\}$ with $\ell \geq 3$,
$f_{h-1} \not= f_h$ and $M(f_{h-1},f_h) \geq 1$ for all $h \in \{1,\dots,\ell\}$,
and $(f_0,f_1)=(f_{\ell-1},f_\ell)=(i,j)$.}
\NEWTEXT{
$A$ is {\em normal} if it is a block diagonal matrix with diagonal blocks $M,D$,
where $D$ is a diagonal matrix with all diagonal entries nonzero, and $M$ has the following two properties:
(i) All diagonal entries are nonzero,
(ii) $M$ is of size $r \geq 2$ and {\em fully indecomposable},
that is, it does not admit an $s \times (r-s)$ zero submatrix, where $s \in \{1,\dots,r-1\}$.}
Finally, $A$ is {\em desirable} if it is tiny or small or special or normal.

Given a matrix $A$ with non-negative integer entries, we associate a graph $G(A)$
on $\{(i,j,t):\, i,j \in \{1,\dots,k\},\,t \in \{1,\dots,A(i,j)\}\}$, where $(i,j,t)$ and $(i',j',t')$ are adjacent if and only if
$i \not= i'$ and $j \not= j'$. 
\NEWTEXT{Figure \ref{F1} shows an example.
\newsavebox{\infigone}
\sbox{\infigone}{\rule[-6mm]{0mm}{14mm}$A=$ {\small $\left(\! \begin{array}{ccc} 1 & 2 & 0 \\ 1 & 3 & 1 \\ 1 & 1 & 1 \end{array} \!\right)$} }
\begin{figure}
  \begin{center} 
    \scalebox{0.4}{\includegraphics{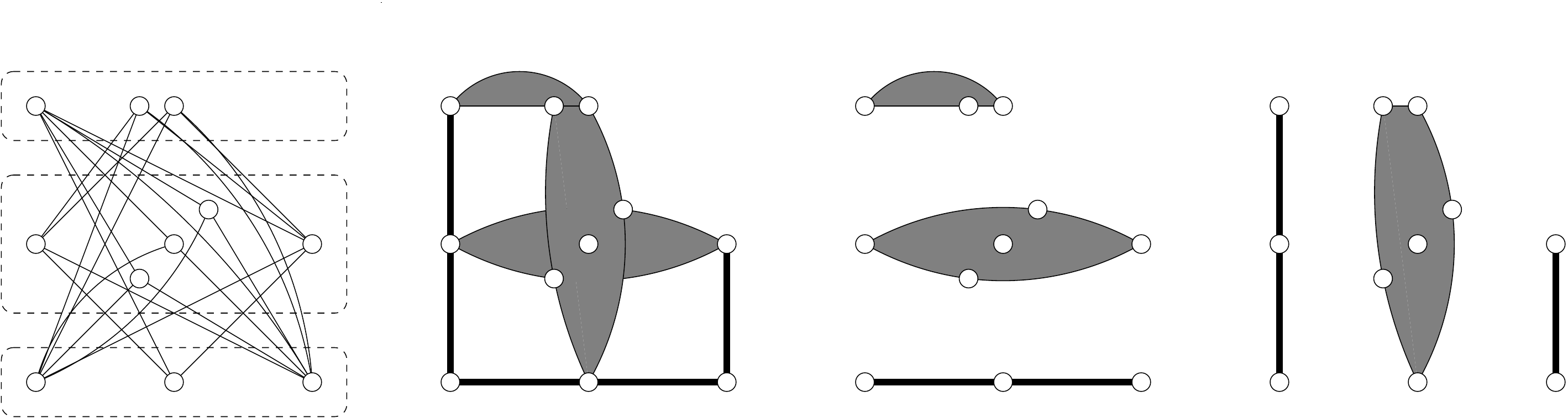}} 
     \caption{\label{F1} The graph $G(A)$ assigned to the fully indecomposable matrix 
     \centerline{\usebox{\infigone}} is drawn leftmost; the dashed boxes indicate one of the two $3$-colorings, the ``horizontal'' one;
     similarly one can identify the ``vertical'' coloring. It may be more convenient to schematically draw the two corresponding overlapping clique partitions of the complementary graph
     $\overline{G}$, as it is done in the second picture from left; the two right figures illustrate the horizontal and the vertical coloring following the same style convention.
     } 
    \end{center}
\end{figure}
}

Our main theorem can be formulated as follows.

\begin{theorem}
  \label{TMaxambkcol}
  Given $k \geq 1$, a graph is maximal ambiguously $k$-colorable if and only if it is isomorphic
  to $G(A)$ for some desirable $k \times k$-matrix $A$.
\end{theorem}

The paper is organized as follows. In the next two sections we give proofs for sufficiency and necessity of the condition
in Theorem \ref{TMaxambkcol} characterizing the maximal ambiguously $k$-colorable graphs.
In Section \ref{STurantype} we derive a {\sc Tur\'an} type result by calculating the maximum number of edges of
an ambiguously $k$-colorable graph on $n$ vertices; moreover, we determine the corresponding extremal graphs (Theorem \ref{TTurantype}).
The situation is somewhat different from {\sc Tur\'an}'s classic Theorem \cite{Turan1941}:
Given $k \geq 3$, there are infinitely many $n$ for which there is only one extremal graph, and infinitely many
$n$ for which there is more than one extremal graph. In Section \ref{SOpenquestions} we discuss generalizations and open questions.

\section{Sufficiency}

In this section, we will prove that whenever $A$ is a desirable $k \times k$-matrix then
$G:=G(A)$ is maximal ambiguously $k$-colorable. According to the definition of being desirable, we distinguish
the cases that $A$ is tiny, small, special, and normal, respectively. %
\NEWTEXT{
For the proof we will need a property of fully indecomposable matrices.

\begin{lemma} \label{sequencelemma}
For every fully indecomposable $r \times r$-matrix $M$ of integers, whenever
$M(i,j) \not=0$ for some $i \not=j$, then there exists a sequence $f_0,\dots,f_\ell$ from $\{1,\dots,r\}$ with $\ell \geq 3$,
$f_{h-1} \not= f_h$ and $M(f_{h-1},f_h) \not= 0$ for all $h \in \{1,\dots,\ell\}$,
and $(f_0,f_1)=(f_{\ell-1},f_\ell)=(i,j)$.
\end{lemma}

{\bf Proof.}
Suppose that $M(i,j) \not= 0$ for some $i \not= j$ and let $S$ be the set of indices $p$ from $\{1,\dots,r\}$ such that
there exists a sequence $f_0,\dots,f_\ell$ from $\{1,\dots,r\}$ with $\ell \geq 1$,
$f_{h-1} \not= f_h$ and $M(f_{h-1},f_h) \not= 0$ for all $h \in \{1,\dots,\ell\}$, $(f_0,f_1)=(i,j)$, and $f_\ell=p$ (*).
Obviously, $j \in S$. If $i \in S$ then there is a sequence as in (*), with $p=i$, and, hence, $\ell \geq 2$, and
the sequence $f_0,f_1,\dots,f_{\ell-1},f_\ell,j$ proves the statement.
So we may assume that $S$ is a nonempty proper subset of $\{1,\dots,r\}$, and so is $T:=\{1,\dots,r\} \setminus S$.
For $p \in S$ and $q \in T$ we infer $M(p,q)=0$, for otherwise there will be a sequence as in (*), and
the sequence $f_0,f_1,\dots,f_{\ell-1},f_\ell,q$ would prove $q \in S$, contradiction.
It follows that $M|S \times T$ is an $|S| \times (r-|S|)$ zero submatrix of $M$,
contradicting the assumption that $M$ is fully indecomposable.
\hspace*{\fill}$\Box$

Let us now turn back to the proof. }%
For any set ${\cal A}$ of pairwise disjoint nonempty sets,
we define {\em the complete ${\cal A}$-partite graph} to be the graph on the vertex $\bigcup {\cal A}$ such that
there is an edge between vertices $a$ and $b$ if and only if $a,b$ are from distinct members of ${\cal A}$.
Any graph isomorphic to the complete ${\cal A}$-partite graph for some ${\cal A}$ with $|{\cal A}|=k$ is called
{\em complete $k$-partite}, and clearly admits a unique $k$-coloring.

{\bf Case 1.} If $A$ is tiny then $G=G(A)$ is a complete graph on at most $k-1$ vertices minus a single edge $xy$. 
It has two distinct $(k-1)$-colorings: One where all vertices form single classes, and another one where
$x,y$ belong to the same class. However, adding the only missing edge $xy$ produces a complete graph
on at most $k-1$ vertices, which is {\em uniquely $k$-colorable}, that is, it has only one $k$-coloring.
Hence $G$ is maximal ambiguously $k$-colorable.

{\bf Case 2.}
If $A$ is small then $G=G(A)$ is complete $(k-1)$-partite, its unique $(k-1)$-coloring ${\cal A}$ consists of sets of size $1$ or $2$,
and at least one member of ${\cal A}$ has size $2$, say $X$.
${\cal A}$ is a $k$-coloring, and by replacing $X$ in ${\cal A}$ by the singletons formed by its two elements
we obtain another $k$-coloring distinct from ${\cal A}$. Hence $G$ is ambiguously $k$-colorable.
Now if $x,y$ are distinct and nonadjacent then $\{x,y\} \in {\cal A}$. 
If ${\cal B}$ is any $k$-coloring of $G+xy$ then we first note that $x,y$ form a clique $K$
of size $k$ together with any selection of vertices $z_A \in A$, $A \in {\cal A}-\{\{x,y\}\}$.
Since every vertex in some $A \in {\cal A}-\{\{x,y\}\}$ is adjacent to all vertices in $K - \{z_A\}$,
it must belong to the same class of ${\cal B}$ as $z_A$. Hence ${\cal B}=({\cal A}-\{\{x,y\}\}) \cup \{\{x\},\{y\}\}$, implying that $G+xy$
is uniquely $k$-colorable. It follows that $G$ is maximal ambiguously $k$-colorable.

{\bf Case 3.}
If $A$ is special then $G=G(A)$ is obtained from a complete $k$-partite graph with unique $k$-coloring ${\cal A}$
by adding a single new vertex $v$ and making it adjacent to all vertices from $\bigcup ({\cal A}-\{S,T\})$,
where we fix $S \not= T$ from ${\cal A}$. It follows that ${\cal A}_S:=({\cal A}-\{S\}) \cup \{S \cup \{v\}\}$ and
${\cal A}_T:=({\cal A}-\{T\}) \cup \{T \cup \{v\}\}$ are distinct $k$-colorings of $G$, so that $G$ is ambiguously $k$-colorable. 
If ${\cal B}$ is any $k$-coloring of $G$ then vertices from distinct classes from ${\cal A}$ (disregarding $v$)
must be in distinct classes from ${\cal B}$,
so that $\{Z-\{v\}:\,Z \in {\cal B}\}$ equals ${\cal A}$.
Therefore, ${\cal B}$ equals either ${\cal A}_S$ or ${\cal A}_T$.
Now suppose that $G'$ is obtained from $G$ by adding a single edge between
two nonadjacent vertices $x,y$. If $x,y$ belong to the same class from ${\cal A}$ then $G'$ has a clique of size $k+1$
and is, therefore, not $k$-colorable. Otherwise, we may assume that $x=v$ and $y \in S$ without loss of generality, and consider
any $k$-coloring ${\cal C}$ of $G'$. Since ${\cal C}$ is a $k$-coloring of $G$, too, it equals either ${\cal A}_S$ or ${\cal A}_T$,
but it cannot be ${\cal A}_S$. Hence ${\cal C}={\cal A}_T$ is the only $k$-coloring of $G'$, proving that $G'$ is
not ambiguously $k$-colorable. Hence $G$ is maximal ambiguously $k$-colorable.

{\bf Case 4.}
Now let $A$ be normal and let $M,D$ be as in the definition of normal.
Without loss of generality, $M=A|\{1,\dots,r\}^2$.
Set $A_i:=\{(i,j,s):\,j \in \{1,\dots,k\},\,s \in \{1,\dots,A(i,j)\}\}$, and
$B_j:=\{(i,j,s):\,i \in \{1,\dots,k\},\,s \in \{1,\dots,A(i,j)\}\}$ for $i,j \in \{1,\dots,k\}$.
Then ${\cal A}:=\{A_1,\dots,A_k\}$ and ${\cal B}:=\{B_1,\dots,B_k\}$ are $k$-colorings of $G:=G(A)$, and $|A_i \cap B_j|=A(i,j)$.
By (i) in the definition of {\em normal}, $|A_i \cap B_i|=A(i,i) \geq 1$ for all $i \in \{1,\dots,k\}$.
From ${\cal A}={\cal B}$ it would follow $A_i=B_i$
for all $i \in \{1,\dots,k\}$ and hence $A(i,j)=|A_i \cap B_j|=0$ for all $i \not= j$ from $\{1,\dots,k\}$ ---
contradicting the fact that $A$ is not a diagonal matrix (since $M$ is fully indecomposable and $r \geq 2$).
This proves ${\cal A} \not= {\cal B}$, hence $G$ is ambiguously $k$-colorable.

Now let ${\cal C}$ be any coloring of $G$. Since vertices from distinct $A_i \cap B_i \not= \emptyset$ must be in distinct classes
of ${\cal C}$, we may list the members of ${\cal C}$ as $C_1,\dots,C_k$ such that $A_i \cap B_i \subseteq C_i$ for all
$i \in \{1,\dots,k\}$. Since every vertex in $A_i \cap B_j$ is adjacent to all vertices in $A_{i'} \cap B_{i'} \not= \emptyset$
for $i' \in \{1,\dots,k\} - \{i,j\}$, we deduce $A_i \cap B_j \subseteq C_i \cup C_j$ for $i,j \in \{1,\dots,k\}$; this statement
strengthens as follows.

{\bf Claim.}
$A_i \cap B_j \subseteq C_i$ or $A_i \cap B_j \subseteq C_j$ for all $i,j \in \{1,\dots,k\}$.

To prove the claim, observe that the statement is obviously true if $|A_i \cap B_j| \leq 1$ or $i=j$. 
Otherwise, $i,j$ are distinct and $|A_i \cap B_j| \geq 2$, so $i,j \in \{1,\dots,r\}$ and $M(i,j) \geq 2$.
\OLDTEXT{We take $\ell$ and $f_0,\dots,f_\ell$ as in (iii) of the definition of {\em normal}.}
\NEWTEXT{We take $\ell$ and $f_0,\dots,f_\ell$ as in Lemma \ref{sequencelemma}.}
Suppose that $A_i \cap B_j \not\subseteq C_i$; then there exists a $y \in A_i \cap B_j \cap C_j$.
We show inductively that $A_{f_{h-1}} \cap B_{f_h}$ is a nonempty subset of $C_{f_h}$ for all $h \in \{2,\dots,\ell\}$.
Let $h \in \{2,\dots,\ell\}$. If $h>2$ then there exists a vertex in $z \in A_{f_{h-2}} \cap B_{f_{h-1}}$
since $M(f_{h-2},f_{h-1}) \geq 1$; by induction, $z \in A_{f_{h-2}} \cap B_{f_{h-1}} \cap C_{f_{h-1}}$.
If, otherwise, $h=2$ then we take $z:=y \in A_{f_0} \cap B_{f_1} \cap C_{f_1}$.
Now consider any vertex $w \in A_{f_{h-1}} \cap B_{f_h} \subseteq C_{f_{h-1}} \cup C_{f_h}$;
there is at least one such vertex, since $M(f_{h-1},f_h) \geq 1$.
Since $f_{h-2}\not=f_{h-1}$ and $f_{h-1} \not=f_h$ we know that $w,z$ are adjacent.
Therefore, they do not belong to the same class from ${\cal C}$, implying that $w \in C_{f_h}$.
This accomplishes the induction.
For $h=\ell$ we get $A_i \cap B_j \subseteq C_j$. This proves the claim.

For $X \in {\cal A}$, set $X^*:=X \times \{\emptyset\}$; for ${\cal X} \subseteq {\cal A}$ define
${\cal X}^*:=\{X^*:\,X \in {\cal X}\}$. 
We construct an auxilary bipartite graph $H$ on ${\cal A}^* \cup {\cal B}$
where $A_i^* \in {\cal A}^*$ and $B_j \in {\cal B}$ are connected by an edge in $H$ if and only if $A(i,j) \geq 1$.
The shape of $A_i,B_j$ ensures that $(x,\emptyset) \not\in V(G)$ for all $x \in V(G)$, so that ${\cal A}^*,{\cal B}$ are disjoint even
if ${\cal A},{\cal B}$ are not.

Let us color an edge $A_i^* B_j$ in $H$ with color $C_i$ if $A_i \cap B_j \subseteq C_i$
and with color $C_j$ if $A_i \cap B_j \subseteq C_j$ (by the claim, every edge of $H$ receives exactly one color).
The set of edges colored with a fixed color $C_\ell$ form a star $H_\ell$ in $G$,
because otherwise there would be disjoint edges $A_i^*B_j$, $A_{i'}^*B_{j'}$ in $H$, colored with $C_\ell$, meaning that
there exists a vertex $v \in A_i \cap B_j$ and a vertex $v' \in A_{i'} \cap B_{j'}$ such that $v,v' \in C_\ell$;
now $v,v'$ are adjacent (as $i\not=i'$ and $j \not= j'$), contradicting the fact that $C_\ell$ is an anticlique.
Moreover, $A_\ell^* B_\ell \in E(H_\ell)$, so that $E(H_1),\dots,E(H_k)$ form a partition of $E(H)$.
Observe that $E(H_\ell)=\{A_\ell^*B_\ell\}$ for $\ell>r$.
For each star $H_\ell$, choose a center $x_\ell$ (which is either $A_\ell^*$ or $B_\ell$).
Let $I:=\{\ell \in \{1,\dots,r\}:\,x_\ell=B_\ell\}$ and let $J:=\{1,\dots,r\}-I$.
It follows that there cannot be an edge $A_i^*B_j$ with $i \in I$ and $j \in J$,
because neither of its endvertices is the center of any star of our star decomposition $H_1,\dots,H_k$.
By definition of $A$ and $M$, we get $M(i,j)=0$ for all $i \in I$ and $j \in J$,
so that $M|I \!\times\! J$ is an $|I| \!\times\! (r-|I|)$ zero submatrix of $M$.
As $M$ is fully indecomposable, $|I| \in \{0,r\}$ follows.

If $|I|=0$ then $x_j=A_j^*$ for all $j \in \{1,\dots,r\}$; hence any edge $A_i^*B_j$ of $H$ with $i,j \in \{1,\dots,r\}$
received color $C_i$, implying that $A_i \cap B_j \subseteq C_i$ for all $i,j \in \{1,\dots,r\}$.
The latter statement extends to {\em all} $i,j \in \{1,\dots,k\}$, since $A_i=B_i$ for $i>r$.
It follows $C_i=A_i$ for all $i \in \{1,\dots,k\}$.
Consequently, ${\cal C}={\cal A}$. --- If, otherwise, $|I|=r$, then ${\cal C}={\cal B}$ follows analogously.

It follows that ${\cal A},{\cal B}$ are the only $k$-colorings of $G$.
Now take any two distinct nonadjacent vertices $x,y$ and suppose that ${\cal C}$ is any $k$-coloring of the graph $G'$
obtained from $G$ by adding a single edge connecting $x,y$. Then ${\cal C}$ is a $k$-coloring of $G$, too,
and hence one of ${\cal A}$ or ${\cal B}$. Since $x,y$ are nonadjacent, they belong either to the same set from ${\cal A}$
or to the same set from ${\cal B}$. In the first case it follows ${\cal C}={\cal B}$ necessarily, and in the second
case we deduce ${\cal C}={\cal A}$; in either case, $G'$ is uniquely $k$-colorable. This proves that
$G$ is maximal ambiguously $k$-colorable.

Hence, $G=G(A)$ is maximal ambiguously $k$-colorable for every desirable $k \times k$-matrix.

\section{Necessity}

Let $G$ be a maximal ambiguously $k$-colorable graph.
Let us prove that $G$ is isomorphic to $G(A)$ for some desirable $k \times k$-matrix $A$.

{\bf Claim 1.}
If ${\cal A}$ is a $(k-1)$-coloring of $G$ then $|A| \leq 2$ for all $A \in {\cal A}$.

Suppose, to the contrary, that $|A| \geq 3$ for some $A \in {\cal A}$, and let $x,y$ be distinct vertices from $A$.
Then $({\cal A}-\{A\}) \cup \{\{x\},A-\{x\}\}$ and $({\cal A} -\{A\}) \cup \{\{y\},A-\{y\}\}$ are distinct $k$-colorings of $G+xy$, contradiction.
This proves Claim 1.

{\bf Claim 2.} If ${\cal A}$ is a $(k-1)$-coloring of $G$ and $|A| \geq 2$ for some $A \in {\cal A}$ then $G$ is complete ${\cal A}$-partite.

By Claim 1, all members of ${\cal A}$ have at most two vertices. We suppose that some member of ${\cal A}$ consists of exactly
two vertices, say, $\{u,v\}$. If there were nonadjacent vertices $x,y$ from distinct classes of ${\cal A}$ then ${\cal A}$ and
$({\cal A}-\{\{u,v\}\}) \cup \{\{u\},\{v\}\}$ were distinct $k$-colorings of $G+xy$, contradiction.
Hence $G$ is complete ${\cal A}$-partite, which proves Claim 2.

{\bf Claim 3.}
If $G$ is $(k-1)$-colorable then there exists a $(k-1)$-coloring ${\cal A}$ such that $G$ is complete ${\cal A}$-partite.

Let ${\cal A}$ be a $(k-1)$-coloring of $G$. By Claim 2 we may assume that all classes of ${\cal A}$ are singletons, so that
$|G|=|{\cal A}| \leq k-1$. We are done if $G$ is complete. So suppose that $x,y$ are nonadjacent vertices. Then
$({\cal A}-\{\{x\},\{y\}\}) \cup \{\{x,y\}\}$ is a $(k-1)$-coloring of $G$, and Claim 2 applies to the modified coloring.
This proves Claim 3.

Now suppose first that $G$ is $q$-colorable for some $q \leq k-2$.
By Claim 3, there exists a $(k-1)$-coloring ${\cal A}$ such that $G$ is complete ${\cal A}$-partite;
in fact, ${\cal A}$ must be a $q$-coloring, for otherwise $G$ would contain $K_{q+1}$, contradicting $q$-colorability.
By Claim 1, the members of ${\cal A}$ have at most two vertices. Suppose, to the contrary, that two distinct members of ${\cal A}$
consist of exactly two vertices each, say, $\{x,y\}$ and $\{u,v\}$. Then $({\cal A}-\{\{x,y\}\}) \cup \{\{x\},\{y\}\}$ and
$({\cal A}-\{\{x,y\},\{u,v\}\} \cup \{\{x\},\{y\},\{u\},\{v\}\}$ are distinct $k$-colorings of $G+xy$ (as $q \leq k-2$), contradiction.
Hence ${\cal A}$ is a $q$-coloring such that $G$ is complete ${\cal A}$-partite, where at most one class has more than one element.
Since $G$ is not complete (as it is not uniquely $k$-colorable), exactly one class has more than one element.
This implies that $G$ is isomorphic to $G(A)$ for some tiny matrix $A$.

Secondly, suppose that $G$ is $(k-1)$-colorable and not $(k-2)$-colorable.
By Claim 3, there exists a $(k-1)$-coloring ${\cal A}$ such that $G$ is complete ${\cal A}$-partite.
By Claim 1, all classes of ${\cal A}$ have at most two vertices, and, as $G$ is not complete as it is ambiguously $k$-colorable,
at least one class of ${\cal A}$ must have exactly two vertices. Since $G$ is not $(k-2)$-colorable, $|{\cal A}|=k-1$. 
But this implies that $G$ is isomorphic to $G(A)$ for some small matrix $A$.

Finally, let us assume that $G$ is not $(k-1)$-colorable, and consider two distinct $k$-colorings ${\cal A}$, ${\cal B}$.
For $X \in {\cal A}$, set $X^*:=X \times \{\emptyset\}$; for ${\cal X} \subseteq {\cal A}$ define
${\cal X}^*:=\{X^*:\,X \in {\cal X}\}$, just as above. Here we may assume without loss of generality that $(x,\emptyset) \not\in V(G)$ for all $x \in V(G)$,
so that ${\cal A}^*,{\cal B}$ are disjoint.
Let $H$ be the bipartite auxilary graph with classes on ${\cal A}^* \cup {\cal B}$,
where there is an edge connecting $A^* \in {\cal A}^*$ and $B \in {\cal B}$
if and only if $A \cap B \not= \emptyset$. If there was an ${\cal X} \subseteq {\cal A}$
such that $|{\cal Y}:=N_H({\cal X}^*)|<|{\cal X}|$
then $\bigcup {\cal X} \subseteq \bigcup {\cal Y}$, and 
$({\cal A}-{\cal X}) \cup \{Y \cap (\bigcup {\cal X}):\, Y \in {\cal Y},\, Y \cap (\bigcup {\cal X}) \not=\emptyset\}$
is a $(k-1)$-coloring of $G$, contradiction. Hence $H$ satisfies the {\sc Hall} condition and, thus, has a perfect matching $M$.
We may label the members of ${\cal A}$, ${\cal B}$ by $A_1,\dots,A_k$ and $B_1,\dots,B_k$ such
that $M=\{A_1^*B_1,\dots,A_k^*B_k\}$ and such that,
for some $r \in \{1,\dots,r\}$, $A_j=B_j$ if and only if $j>r$.
Since $A_j=B_j$ for all $j > 1$ would imply $A_1=B_1$ (contradicting ${\cal A} \not= {\cal B}$), we know $r \geq 2$.

Let us define a $k \times k$ matrix $A$ with nonnegative integer entries by $A(i,j):=|A_i \cap B_j|$ for $i,j \in \{1,\dots,k\}$. 
It is obvious that $G$ is isomorphic to $G(A)$, and that all diagonal entries are at least $1$.
We claim that $A$ is either special or normal.

Suppose first that $A_i \subseteq B_i$ holds for some $i \in \{1,\dots,r\}$.
Since $A_i \not= B_i$, there exists a vertex $x$ in $B_i-A_i$; clearly, $x \in B_i \cap A_j$ for some $j \not= i$,
$x$ is not adjacent to any of $A_i$, and $A_j$ contains a vertex distinct from $x$, say $y \in A_j \cap B_j$.
Consequently, ${\cal C}:=({\cal A}-\{A_i,A_j\}) \cup \{A_i \cup \{x\},A_j-\{x\}\}$ is a $k$-coloring of $G$.
Observe that $x,y$ are non-adjacent and both ${\cal C}$ and ${\cal B}$ are colorings of the graph $G'$ obtained from $G$
by adding a single edge connecting $x$ and $y$. Since $G$ is maximal ambiguously $k$-colorable,
${\cal C}={\cal B}$ follows, implying that $A$ is special, where $A(j,i)=1$ is the unique nonzero off-diagonal entry. 
Analogously, if $B_i \subseteq A_i$ holds for some $i \in \{1,\dots,r\}$ then $A$ is special, too.

Hence, for the remaining argument, we may assume that $A_1,\dots,A_r$, $B_1,\dots,B_r$ are incomparable with respect to $\subseteq$.
Set $M:=A|\{1,\dots,r\}^2$ and $D:=A|\{r+1,\dots,k\}^2$. We claim that $A$ is normal.
\OLDTEXT{Since $A$ has no zero diagonal entries, it suffices to verify conditions (ii) and (iii) to $M$ in the definition of {\em normal}.}
\NEWTEXT{Since $A$ has no zero diagonal entries, it suffices to verify condition (ii) to $M$ in the definition of {\em normal}.}

Suppose, to the contrary, that $M$ is not fully indecomposable. Then there exist nonempty $I,J \subseteq \{1,\dots,r\}$
such that $|I|+|J|=r$ and $M|I \times J$ is zero everywhere, i. e. $A_i^*,B_j$ are not connected by an edge in $H$
whenever $i \in I$ and $j \in J$. Since $A_i^*B_i \in E(H)$ for all $i \in \{1,\dots,r\}$,
$I$ and $J$ are disjoint and, therefore, form a partition of $\{1,\dots,r\}$.
Setting ${\cal X}:=\{A_i:\, i \in I\}$ and ${\cal Y}:=N_H({\cal X}^*) \subseteq \{B_j:\, j \in I\}$,
we see again that $\bigcup {\cal X} \subseteq \bigcup {\cal Y}$, and 
${\cal D}$:=$({\cal A}-{\cal X}) \cup \{Y \cap (\bigcup {\cal X}):\, Y \in {\cal Y},\, Y \cap (\bigcup {\cal X}) \not=\emptyset\}$
is a $k$-coloring of $G$. 
\NEWTEXT{Figure \ref{F2} illustrates the construction; $A_1,\dots,A_r$ constitutes the ``horizontal'' coloring, $B_1,\dots,B_r$ the ``vertical'' one.
\newsavebox{\infigtwo}
\sbox{\infigtwo}{\rule[-11mm]{0mm}{24mm}$A=$ {\small $\left(\! \begin{array}{ccccc} 1 & 1 & 1 & 1 & 1 \\ 1 & 1 & 1 & 1 & 1 \\ 0 & 0 & 1 & 1 & 1 \\ 0 & 0 & 1 & 1 & 1 \\ 0 & 0 & 1 & 1 & 1 \end{array} \!\right)$} }
\begin{figure}
  \begin{center} 
    \scalebox{0.4}{\includegraphics{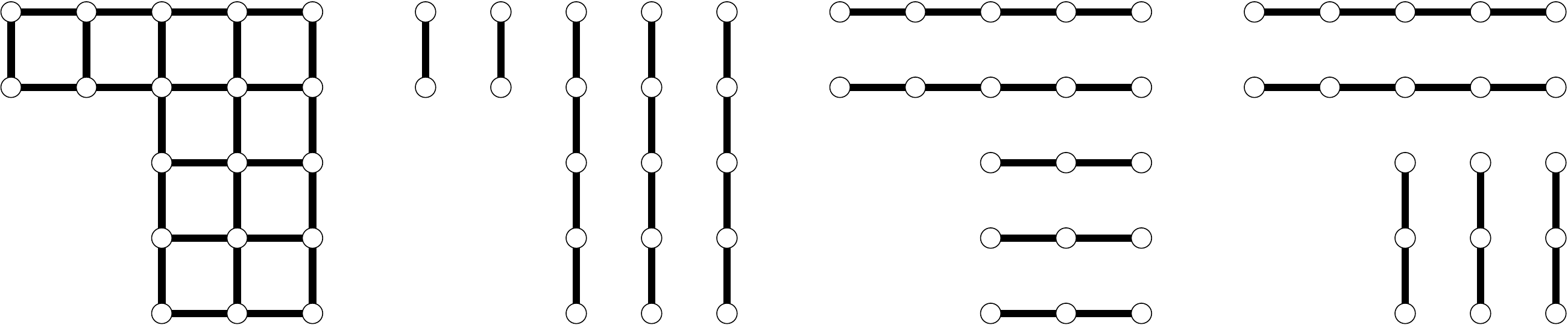}} 
     \caption{\label{F2}
       How to obtain a third coloring (rightmost) from the the vertical and horizontal colorings (middle) of the graph $G(A)$ assigned to the matrix 
       \centerline{\usebox{\infigtwo},} \medskip which is not fully indecomposable.
       The graph (leftmost) and the colorings are schematically drawn as clique covers of $\overline{G(A)}$. }
    \end{center}
\end{figure}}
Take $A_i \in {\cal X}$. Since $A_i$ is not contained in $B_i \in {\cal Y}$, 
there must be a $B_j \in {\cal Y}$ with $j \not= i$ and distinct nonadjacent vertices $x \in A_i \cap B_i$, $y \in A_i \cap B_j$.
Now ${\cal B}$ and ${\cal D}$ are colorings of the graph $G'$ obtained from $G$ by adding a single edge connecting $x,y$.
Take any $i' \in J \not=\emptyset$; then $A_{i'} \in {\cal A}-{\cal X} \subseteq {\cal D}$.
Since $B_{i'}$ is not contained in $A_{i'}$ 
we see that ${\cal B}$ is distinct from ${\cal D}$, contradicting the fact that $G$ is maximal ambiguously $k$-colorable.
This proves that $M$ is fully indecomposable.

This accomplishes the proof of Theorem \ref{TMaxambkcol}.

\section{A Tur\'an type consequence}
\label{STurantype}

Given integers $r,n$, the {\em Tur\'an number} of $n$ and $K_{r+1}$ is the largest number $\ex(n,K_{r+1})$
of edges a (simple) graph on $n$ vertices without $K_{r+1}$ as a subgraph can have,
and graphs on $n$ vertices without a $K_{r+1}$ as a subgraph and with $\ex(n,K_{r+1})$ edges
are called {\em $(n,K_{r+1})$-extremal}.
For $n \leq r$, $\ex(n,K_{r+1})={n \choose 2}$, and $K_n$ is the only extremal graph up to isomorphism,
whereas, for $n>r$, the only extremal graph up to isomorphism is the balanced complete $r$-partite graph $T(n,r)$ on $n$ vertices,
that is: $T(n,r)$ is complete ${\cal A}$-partite where $|{\cal A}|=r$,
$|\bigcup {\cal A}|=n$, and $||A|-|A'|| \leq 1$ for all $A,A'$ from ${\cal A}$
(and $\ex(n,K_{r+1})=|E(T(n,r))|$, which has various algebraic representations \cite{Turan1941}.

Let us call a graph {\em maximal $k$-colorable} if it is $k$-colorable but any graph obtained from $G$
by adding a single edge between two distinct nonadjacent vertices is not. 
Obviously, the maximal $k$-colorable graphs on $n$ vertices are complete if $n \leq k$,
and complete ${\cal A}$-partite for some $k$-coloring ${\cal A}$ with $|{\cal A}|=k$
if $n \geq k$, so that $\ex(n,K_{k+1})$ is equal to
the largest number of edges a $k$-colorable graph on $n$ vertices can have.

Let us determine the corresponding extremal numbers for the property of being
ambiguously $k$-colorable. We start with the following Lemma.

\begin{lemma}
  \label{LTurantype}
  Let ${\cal A}$ be a partition of order $k$ of a set of order $n$, and let $G$ be a spanning subgraph of
  the complete ${\cal A}$-partite graph. Let $\alpha:=\lfloor \frac{n}{k} \rfloor$ and suppose that
  $A_1,\dots,A_r$ are members of ${\cal A}$ of order at most $\alpha$. Fix $H:=G(\bigcup_{j=1}^r A_j)$,
  and let $d$ be the number of those edges of the complete $\{A_1,\dots,A_r\}$-partite graph which are not in $E(H)$.
  Let $r_0$ be the number of sets among $A_1,\dots,A_r$ with at most $\alpha-1$ elements.
  Then $|E(G)| \leq {\mathrm ex}(n,K_{k+1})-( 2 \cdot (\alpha \cdot r - |V(H)|) - r_0)-d$.
\end{lemma}

{\bf Proof.}
Let $A_1,\dots,A_k$ be the members of ${\cal A}$.
As long as there exists a $j \in \{1,\dots,r\}$ such that $|A_j| \leq \alpha-1$,
we modify $G$ and the partition ${\cal A}$ --- but not $H$ --- such that each step
preserves $G(V(H))$, $V(H) \cap A_1,\dots,V(H) \cap A_r$,
$|A_i| \leq \alpha$ for all $i \in \{1,\dots,r\}$,
and the size of $G$ increases:
First observe that there exists an $i \in \{1,\dots,k\}$ such that $|A_i| \geq \lceil \frac{n}{k} \rceil \geq \lfloor \frac{n}{k} \rfloor =\alpha$.
If $|A_i| \leq \alpha$ for all $i \in \{1,\dots,k\}$ then $k$ divides $n$ and $|A_i|=\frac{n}{k}$ for all $i \in \{1,\dots,k\}$,
contradicting the existence of $A_j$ as above. Hence $|A_i| \geq \alpha+1$ holds for some $i \in \{1,\dots,k\}$.
Clearly, $i \not\in \{1,\dots,r\}$.
Choose $x \in A_i$, delete all edges of $E_G(\{x\},A_j)$ from $G$, add a single edge from $x$ to each $y \in A_i-\{x\}$,
and call the resulting graph $G'$. Set $A'_i:=A_i-\{x\}$ and $A'_j:=A_j \cup \{x\}$ and $A'_p:=A_p$ for all $p \in \{1,\dots,k\}-\{i,j\}$.
${\cal A}':=\{A'_1,\dots,A'_k\}$ is a partition of $G'$ into anticliques, $A'_1,\dots,A'_r$ have order at most $\alpha$,
$G'(V(H))=H$, and $V(H) \cap A'_1=V(H) \cap A_1,\dots,V(H) \cap A'_r=V(H) \cap A_r$.
Observe that $|E(G')| \geq |E(G)|-|A_j|+|A_i-\{x\}|$, so that in each step where $|A_j| \leq \alpha-2$ we get
at least two additional edges, whereas in each step where $|A_j|=\alpha-1$ we get at least one. 
This way, we may perform a total of $\alpha \cdot r-|V(H)|$ steps, at most $r_0$ of which increase the order of some $A_j$
from $\alpha-1$ to $\alpha$. By finally adding $d$ edges between those pairs of non-adjacent vertices from $V(H)$
which are in distinct classes from the finally constructed partition, we obtain a graph which is still $k$-partite, but
gain at least $( 2 \cdot (\alpha \cdot r - |V(H)|) - r_0)+d$ in size compared to the initial graph.
\hspace*{\fill}$\Box$

A desirable $k \times k$-matrix $A$ is called {\em row-sum-balanced} if
$|\sum_{j=1}^k A(i,j) - A(i',j)|$ $\leq$ $1$ for all $i,i' \in \{1,\dots,k\}$, that is, the difference of any two row-sums is $0$ or $\pm1$.
Likewise, $A$ is {\em column-sum-balanced} if
$|\sum_{i=1}^k A(i,j) - A(i,j')| \leq 1$ for all $j,j' \in \{1,\dots,k\}$.
$A$ is {\em balanced} if it is both row- and column-sum-balanced.

Let us calculate the number of edges of $G(A)$ for some desirable $k \times k$-matrices $A$, given that 
$n:=|G(A)|$ (which is the sum over all entries of $A$).

If $A$ is tiny then $n \leq k-1$, and $|E(G(A))|={n \choose 2}-1$ $=$ $\ex(n,K_{k+1})-1$.
If $A$ is small then $G(A)$ is a complete $(k-1)$-partite graph, and we get $k \leq n \leq 2k-2$;
it follows $|E(G(A))|={n \choose 2}-(n-k+1)=\ex(n,K_{k+1})-1$.

Let $A_i:=\{(i,j,\ell):\, j \in \{1,\dots,k\},\ell \in \{1,\dots,A(i,j)\}\}$ and ${\cal A}:=\{A_1,\dots,A_k\}$.
Suppose that $A$ is row-sum-balanced, that is, $||A_i|-|A_j|| \leq 1$ for all $i,j \in \{1,\dots,k\}$.
Equivalently, each $|A_i|$ is one of $\alpha:=\lfloor \frac{n}{k} \rfloor$, $\lceil \frac{n}{k} \rceil$.
If $A$ is special then let $i \not= j$ be the unique indices from $\{1,\dots,k\}$ such that $A(i,j)=1$.
Then $G(A)$ is obtained from the complete $\{A_1,\dots,A_k\}$-partite graph by deleting
all edges connecting $(i,j,1)$ to $(j,j,1),\dots,(j,j,A(j,j))$. 
Therefore, $|E(G(A))|=\ex(n,K_{k+1})-|A_j|$, which is at most $\ex(n,K_{k+1})-\alpha$.
Accordingly, let us call a special matrix $A$ {\em (a)-special} if it is row-sum-balanced and
the sum of the entries in row $j$ is $\lfloor \frac{n}{k} \rfloor$, where
$j$ is the index of the unique column with an off-diagonal entry.
Observe that we can realize an (a)-special matrix for all $n \geq k+1$.
If $k$ divides $n$ or $n \leq 2k-1$ then, up to isomorphism, they induce all one and the same graph,
whereas if $k$ does not divide $n$
and $n \geq 2k$ then, up to isomorphism, they induce all one among two graphs, depending on whether $(i,j,1)$ is non-adjacent to
$\lfloor \frac{n}{k} \rfloor+\lfloor \frac{n}{k} \rfloor-1$ or to $\lfloor \frac{n}{k} \rfloor + \lceil \frac{n}{k} \rceil -1$ many other vertices ---
corresponding to the cases that $|A_i|=\alpha$ and $|A_i|=\alpha+1$, respectively.
For the special case that $k+1 \leq n \leq 2k-1$ we get $|E(G(A))|={n \choose 2}-(n-k)-1$ $=$ $\ex(n,K_{k+1})-1$,
which is equal to the size of a graph on $n$ vertices induced by a small $k \times k$-matrix when $k+1 \leq n \leq 2k-2$.

Symmetrically, let us call a special matrix $A$ {\em (b)-special} if it is colum-sum-balanced and
the sum of the entries in column $j$ is $\lfloor \frac{n}{k} \rfloor$, where
$j$ is the index of the unique column with an off-diagonal entry.
Obviously, $G(A)$ is isomorphic to $G(A^\top)$, so that we got no further graphs this way.
However, there is a third way to realize an ambiguously $k$-colorable graph on $n$ vertices with
$\ex(n,K_{k+1})-\lfloor \frac{n}{k} \rfloor$ edges by a special matrix $A$ which is neither row- nor column-sum-balanced:
Let us call a special matrix $A$ {\em (c)-special}, if $A(i,i)=A(j,j)=\alpha -1$ where $i,j$ are the unique
indices $i \not= j$ with $A(i,j) \not=0$, $A(\ell,\ell) \in \{\alpha,\alpha+1\}$ for all $\ell$ distinct from $i,j$,
and $A(\ell,\ell)=\alpha+1$ for at least one $\ell$. 
Let $A'$ be obtained from $A$ by adding $1$ at position $(j,j)$ and subtracting $1$ at position $(\ell,\ell)$, where
$A(\ell,\ell)=\alpha+1$. The row-sums of $A'$ are $\alpha$ or $\alpha+1$, so that $A'$ is (a)-special;
in fact, since the row-sums of at least three rows are $\alpha$, we see that, necessarily,
$n \geq 2k \geq 6$ and $n$ is not congruent $-1$ or $-2$ modulo $n$, and this is also sufficient for the existence of
a (c)-special $k \times k$-matrix whose entries sum up to $n$.
Let us compare the sizes of $G(A)$ and $G(A')$:
We could think of $G(A')$ as obtained from $G(A)$ by deleting vertex $(\ell,\ell,\alpha+1)$ and
adding a new vertex $(j,j,\alpha)$ and connect according to the rules defining $G(A')$;
by the deletion we loose $\alpha$ edges in the complementary graph, whereas by
the addition we gain $\alpha-1$ edges from $(j,j,\alpha)$ to $(j,j,\beta)$, $\beta<\alpha$
plus one further edge from $(j,j,\alpha)$ to $(i,j,1)$ (again in the complementary graph).
Therefore, $|E(G(A'))|=|E(G(A))|$.
Let us summarize by defining $A$ to be {\em very special}, if it is (a)-, (b), or (c)-special.

Suppose now that $A$ is a normal matrix, and let $M,D$ be as in the definition of {\em normal}.
We call $A$ {\em mininormal}, if $A$ is balanced, $M={1 \, 1 \choose 1 \, 1}$, and $2k \leq n < 3k$.
Consequently, the diagonal entries of $D$ are either $2$ or $3$, so that $n \leq 3k-2$.
Up to permutation of the diagonal entries of $D$, there is exactly one mininormal matrix for each $n \in \{2k,\dots,3k-2\}$,
and $G(A)$ is obtained from the complete ${\cal A}$-partite graph $H$ (${\cal A}$ as above) by deleting the two edges of a $1$-factor
of $G(A_1 \cup A_2)$, where $A_1 \not= A_2$ are from ${\cal A}$ of order $2$.
Consequently, $|E(G(A))|={\mathrm ex}(n,K_{k+1})-2$, which is equal to the size of a graph 
on $n$ vertices induced by a very special $k \times k$-matrix.

It turns out that these constructions produce exactly the ambiguously $k$-colorab\-le graphs with the largest number of edges,
as stated in the following theorem.

\begin{theorem}
  \label{TTurantype}
  Let $n,k \geq 2$ be integers.
  Then the maximum number of edges in an ambiguously $k$-colorable graph on $n$ vertices is
  $\ex(n,K_{k+1})-\max\{1,\lfloor \frac{n}{k} \rfloor\}$.    
  The graphs where the bound is attained are isomorphic to $G(A)$, where $A$ is a desirable $k \times k$-matrix such that
  $A$ is tiny or small or very special or mininormal.
\end{theorem}

{\bf Proof.} 
Let $f(n,k):=\max\{|E(G)|:\,G$ is an ambiguously $k$-colorable graph on $n$ vertices$\}$,
and let $g(n,k):=\ex(n,K_{k+1})-\max\{1,\lfloor \frac{n}{k} \rfloor\}$. We have seen before that
$f(n,k) \geq g(n,k)$. Set $\alpha:=\lfloor \frac{n}{k} \rfloor$.

Suppose that $G$ is an ambiguously $k$-colorable graph on $n$ vertices with $f(n,k)$ edges.
Then $G$ is maximal ambiguously $k$-colorable. 
By Theorem \ref{TMaxambkcol}, we may assume that $G=G(A)$ for some desirable matrix $A$.

It remains to show that $|E(G)| < g(n,k)$ or $A$ is tiny, small, very special or mininormal (in these cases,
we know from the above considerations that $|E(G)|=g(n,k)$). Thus, it suffices to analyze the cases
that $A$ is special or normal, respectively.

Suppose first that $A$ is special.
Let $A_i:=\{(i,j,\ell):\, j \in \{1,\dots,k\},\ell \in \{1,\dots,A(i,j)\}\}$ and ${\cal A}:=\{A_1,\dots,A_k\}$.
Let $i,j$ be the unique indices with $i \not= j$ and $A(i,j) \not=0$, let $v:=(i,j,1)$, let $S:=A_i$, and let $T:=A_j$.
Then $G$ is obtained from the complete ${\cal A}$-partite graph $H$ 
by deleting all edges of $E_H(\{v\},T)$ in $H$. If $X$, $Y$ are distinct classes from ${\cal A}$ such that
$\Delta:=|Y|-|X| \geq 2$ then choose any $Z \subseteq Y-\{v\}$ with $|Z|=\lfloor \frac{\Delta}{2} \rfloor$
and define $X':=X \cup Z$, $Y':=Y-Z$, $C':=C$ for all $C \in {\cal A}-\{X,Y\}$;
set ${\cal A}':=\{C':\,C \in {\cal A}\}$,
and obtain $G'$ from the complete ${\cal A}'$-partitite graph $H'$ by deleting all edges from $E_{H'}(\{v\},T')$. 
Observe that $G'$ is isomorphic to a graph $G(A')$ for another special matrix $A'$.

Let us compare the sizes of $G$ and $G'$:
$|E(G')|=|E(H')|-|T'|=|E(H)|+\lfloor \frac{\Delta}{2} \rfloor \cdot \lceil \frac{\Delta}{2} \rceil-|T'|$
$=$ $|E(G)|+ |T| + \lfloor \frac{\Delta}{2} \rfloor \cdot \lceil \frac{\Delta}{2} \rceil-|T'|$
$\geq$ $|E(G)| + \lfloor \frac{\Delta}{2} \rfloor \cdot \lceil \frac{\Delta}{2} \rceil-|Z|$
$=$ $|E(G)| + \lfloor \frac{\Delta}{2} \rfloor \cdot (\lceil \frac{\Delta}{2} \rceil-1)$
$\geq$ $|E(G)|$, 
where the first inequality is an equality if and only if $X=T$,
and the second one is an equality if and only if $\Delta=2$.
Since $f(n,k) \geq |E(G')| \geq |E(G)|=f(n,k)$, we deduce $X=T$ and $\Delta=2$.

It follows that $||X|-|Y|| \leq 1$ and $|Y|-|T| \in \{-1,0,1,2\}$ for all $X,Y \in {\cal A}-\{T\}$.

If $A$ is row-sum-balanced then $|X|-|T| \in \{-1,0,1\}$ for all $X \in {\cal A}-\{T\}$, too,
$H$ is $(n,K_{k+1})$-extremal, and $|E(G)|=|E(H)|-|T| \leq E(H)-\lfloor \frac{n}{2} \rfloor$ $=$ $g(n,k)$, where equality holds
if and only if $|T|=\lfloor \frac{n}{2} \rfloor$; in that case, $A$ is (a)-special.
Analogously, if $A$ is column-sum-balanced, it must be (b)-special.
Hence we may assume that $A$ is neither row- nor column-sum-balanced.
Recall that $T=A_j$. If the sum over each row of $A$ distinct from $j$ would equal the same value $\beta$,
then the sum over row $j$ equals $\beta-2$, and the sum over column $i,j$ equals $\beta-1$,
respectively, whereas the sum over each other column equals $\beta$;
hence $A$ is column-sum-balanced, contradiction.
Otherwise, the sums over the rows distinct from row $j$ take values, say, $\beta$ and $\beta+1$.
The sum over row $j$ is $\beta-1$.
Suppose that $r$ is the number of rows whose sum is $\beta+1$. Then $0<r<k-1$.
Therefore, we get $n=(\beta-1)+r \cdot (\beta+1)+(k-1-r) \cdot \beta$ $=$ $k \cdot \beta + (r-1)$,
implying $\beta=\lfloor \frac{n}{k} \rfloor = \alpha$. If the sum over row $i$ is $\beta+1$, then the sum over each column
equals $\beta$ or $\beta+1$, and we are done. Hence the sum over row $i$ is $\beta=\alpha$, which implies
$A(i,i)=\alpha-1$. Analogously, $A(j,j)=\alpha-1$, and it straightforward to check that $A$ is (c)-special.

Finally, suppose that $A$ is normal and let $M,D$ be matrices as in the definition of {\em normal},
where $M=A|\{1,\dots,r\}^2$ without loss of generality.
Let $A_i,B_j$ and ${\cal A},{\cal B}$ be as in the sufficiency proof of Theorem \ref{TMaxambkcol},
and let $H$ be the complete ${\cal A}$-partite graph.
Observe that $H$ is a supergraph of $G$ and $|E(H)| \leq \mbox{ex}(n,K_{k+1})$.

Consider $j \in \{1,\dots,r\}$, and let ${\cal C}_j:=\{A_i \cap B_j:\, i \in \{1,\dots,r\},\,M(i,j) \geq 1\}$.
${\cal C}_j$ is a partition of $B_j$.
Since $M$ is fully indecomposable, $|{\cal C}_j| \geq 2$.
Consequently, the complete ${\cal C}_j$-partite graph $H_j$ is connected.
Since the edges of every $H_j$ are present in $H$ but not in $G$, 
and since $|E(H)|-|E(G)| \leq \ex(n,K_{k+1})-f(n,k) \leq \max\{1,\alpha\}=\alpha$ (as $n \geq k$),
we may assume that,
for any selection $K \subseteq \{1,\dots,r\}$,
$|K| \leq \sum_{j \in K} (|B_j|-1) \leq \sum_{j \in K} |E(H_j)| \leq \alpha$. 
For $|K|=1$ we obtain $2 \leq |B_j| \leq \alpha+1$ for every $j \in \{1,\dots,r\}$, and by extending
to a set $K:=\{j,j'\}$ of order $2$ (note that $r \geq 2$) we deduce $|B_j|+|B_{j'}| \leq \alpha+2$, so that, in fact,
$|B_j| \leq \alpha$. Taking $K=\{1,\dots,r\}$ and setting $d:=\sum_{j=1}^r |E(H_j)|$,
we get $2 \cdot  r \leq \sum_{i=1}^r |A_i| = \sum_{j=1}^r |B_j| \leq \sum_{j=1}^r (|E(H_j)|+1) = d+r \leq \alpha+r$.
In particular, $r \leq \alpha$. Since the arguments apply symmetrically to the situation where
the roles of the two $k$-colorings are swapped, we get $|A_j| \leq \alpha$ for all $j \in \{1,\dots,r\}$.

Observe that $d$ is the number of edges in the complete $\{A_1,\dots,A_r\}$-partite graph which are not in $G(\bigcup_{i=1}^r A_i)$.
Now Lemma \ref{LTurantype} applies, yielding $|E(G)| \leq {\mathrm ex}(n,K_{k+1})-( 2 \cdot (\alpha \cdot r - m) - r_0)-d$,
where $m:=\sum_{i=1}^r |A_i|$ and $r_0$ is the number of $A_1,\dots,A_r$ of order at most $\alpha-1$. 
It thus suffices to show that
$( 2 \cdot (\alpha \cdot r - m) - r_0)+d \geq \alpha$ and to analyze the case when equality holds.
We estimate
$( 2 \cdot (\alpha \cdot r - (d+r)) - r_0)+d$
$\geq$ $( 2 \cdot (\alpha \cdot r - (d+r)) - r)+d$
$=$ $2 \cdot \alpha \cdot r - d -3 \cdot r$
$\geq$ $2 \cdot \alpha \cdot r -\alpha -3 \cdot r$.
The latter ist at least $\alpha$ if and only if
$2 \cdot \alpha \cdot r - 2 \alpha -3 \cdot r \geq 0$, which is equivalent to
$(\alpha-\frac{3}{2}) \cdot (r -1) \geq \frac{3}{2}$. For $r \geq 3$ (and hence $\alpha \geq r \geq 3$),
we get even the strict inequality $(\alpha-\frac{3}{2}) \cdot (r -1) \geq \frac{3}{2} \cdot 2 > \frac{3}{2}$
(strictness propagates back yielding $|E(G)|<g(n,k)$).

We thus may assume that $r=2$. For $\alpha \geq 4$, we get, again, 
the strict inequality $(\alpha-\frac{3}{2}) \cdot (r -1) \geq \frac{5}{2} \cdot 1 > \frac{3}{2}$,
and hence $|E(G)|>f(n,k)$.
So let us assume $\alpha \leq 3$.
Since $M$ is a fully indecomposable $2 \times 2$-matrix, all entries are nonzero.
Observe that the sum of all entries of $M$ equals $m=|A_1|+|A_2| \leq d+r \leq \alpha+2 \leq 5$,
and $|A_1|,|A_2| \geq 2$, hence $\alpha \in \{2,3\}$.

For $\alpha=3$, we get the case that both $A_1,A_2$ have order $2$
and the case that one has order $2$ and the other has order $3$.
In the first case, $m=4$ and $r_0=2$ and $d=2$, so that
$( 2 \cdot (\alpha \cdot r - m) - r_0)+d = (2 \cdot (3 \cdot 2 - 4) - 2)+2=4>3=\alpha$;
otherwise, $m=5$ and $r_0=1$ and $d=3$, so that
$( 2 \cdot (\alpha \cdot r - m) - r_0)+d = (2 \cdot (3 \cdot 2 - 5) - 1)+3=4>3=\alpha$.
In either case, $|E(G)|<g(n,k)$.

For $\alpha=2$ observe that $|A_1|=|A_2|=2$, and $2k \leq n \leq 3k-1$.
If the graph obtained from $G$ by adding the two missing edges between $A_1$ and $A_2$
was not $(n,K_{r+1})$-extremal, then $|E(G)| < g(n,k) \leq f(n,k)$. Hence
$G$ is obtained from an $(n,K_{r+1})$-extremal graph by deleting a matching between
two classes of size $2$, so that $A$ is mininormal.
\hspace*{\fill}$\Box$

\section{Open questions}
\label{SOpenquestions}

Theorem \ref{TMaxambkcol} implies that the chromatic number $\chi(G)$ of a maximal ambiguously $k$-colorable graph $G$ is equal to
its clique number $\omega(G)$; in particular these graphs satisfy the statement of {\sc Hadwiger}'s Conjecture (see introduction). 
However, it turns out that much more is true:
Let $A$ be any matrix with nonnegative integer entries and let $G(A)$ as defined in the introduction.
If $A$ is a $k \times k$-matrix where all entries of $A$ are constantly $\ell$ then $G(A)$ is isomorphic to
$G^+:=\overline{K_k \times K_k}[\overline{K_\ell}]$.
Here $\overline{G}:=(V(G),\{xy:\, x \not= y$ in $V(G)$, $xy \not\in E(G)\}$ denotes the {\em complementary graph} of a graph $G$,
$G \times H:=(V(G) \times V(H),\{(x,y)(x',y'):\,x=x' \wedge yy' \in E(H)$ or $y=y' \wedge xx' \in E(G)\}$ denotes the
{\em cartesian product} of two graphs $G$ and $H$, and
$G[H]:=(V(G) \times V(H),\{(x,x')(y,y'):\,xx' \in E(G)$ or $x=x' \wedge yy' \in E(H)\}$ denotes their {\em lexicographic product}. ---
Likewise, if all entries of $A$ are at most $\ell$
then $G(A)$ is an induced subgraph of $G^+$.
By Theorem 2.6.(iv) and Theorem 4.1 of \cite{RavindraParthasarathy1978} plus {\sc Lov\'asz}'s celebrated Perfect Graph Theorem
that complements of perfect graphs are perfect \cite{Lovasz1972}, it follows that $G^+$ is perfect, and this is inherited to $G(A)$.
In particular we get the following:

\begin{theorem}
  Every maximal ambiguously $k$-colorable graph is perfect.
\end{theorem}

As sketched in the introduction, let us call a graph {\em $d$-fold $k$-colorable} if it has $d$ pairwise distinct $k$-colorings,
and {\em maximal $d$-fold $k$-colorable} if it is $d$-fold $k$-colorable but adding any edge
between distinct nonadjacent vertices produces a graph which is not.
For a mapping $A:\{1,\dots,k\}^d \rightarrow \mathbb{Z}_{\geq 0}$ let us define a graph $G(A)$
on $\{(i_1,\dots,i_d,s):\,i_j \in \{1,\dots,k\} \mbox{ for }j \in \{1,\dots,d\},\,s \in \{1,\dots,A(i_1,\dots,i_d)\}\}$ where there is
an edge connecting $(i_1,\dots,i_d,s)$ and $(i'_1,$ $\dots,$ $i'_d,s')$ if and only if $i_j \not= i'_j$ for all $j \in \{1,\dots,d\}$.
It is then easy to see that every maximal $d$-fold $k$-colorable graph is isomorphic to a graph $G(A)$
for a suitable mapping $A$ from $\{1,\dots,k\}^d$ to the nonnegative integers.
Now every such graph $G(A)$ is an induced subgraph of the graph
\[ G^+:=\overline{\underbrace{K_k \times \dots \times K_k}_{\mbox{$d$ times}}}[\overline{K_\ell}],\]
where $\ell:=\max\{A(i_1,\dots,a_d):\,i_j \in \{1,\dots,k\} \mbox{ for }j \in \{1,\dots,d\}\}$.
However, $G^+$ is {\em not} perfect for $d>2$, as it has been proved in \cite{RavindraParthasarathy1978}.

A number of natural problems arise. The most difficult ones are perhaps the generalizations of
Theorem \ref{TMaxambkcol} and Theorem \ref{TTurantype}:

\begin{problem}
  Characterize the mappings $A:\{1,\dots,k\}^d \rightarrow \mathbb{Z}_{\geq 0}$ for which $G(A)$ is maximal $d$-fold $k$-colorable.
\end{problem}

\begin{problem}
  Determine the maximum number of edges in a maximal $d$-fold $k$-colorable graph on $n$ vertices.
  Determine the graphs attaining this maximum.
\end{problem}

More particular one could ask if the extremal graphs in this problem are necessarily ``slightly perturbed {\sc Tur\'an}-graphs''
as in the $2$-dimensional case (Theorem \ref{TTurantype}).



We have seen that maximal $d$-fold $k$-colorable graphs are not necessarily perfect for $d>2$. Still one could think that
a maximal $d$-fold $k$-colorable graph has a clique of order $\chi(G)$, but this is not true in general; the argument is based on
the fact that if a triangle free graph $G$ has a perfect matching then $\chi(\overline{G})=|V(G)|/2$ and the perfect matchings
of $G$ correspond to the $\chi(\overline{G})$-colorings of $\overline{G}$. Now take $G$ to be the graph obtained from $K_4$ by subdividing
each of the two edges of some fixed matching twice
\NEWTEXT{(see Figure \ref{F3});
\begin{figure}
  \begin{center} 
    \scalebox{0.4}{\includegraphics{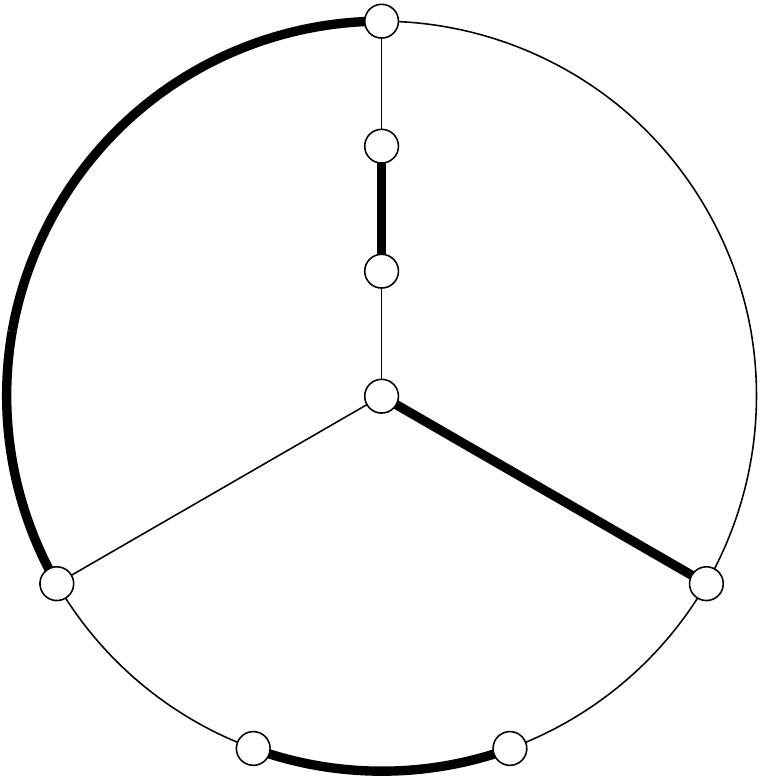}} 
     \caption{\label{F3}  An example of a triangle free graph with exactly $3$ perfect matchings (one is displayed fat); deleting any edge will kill at least one of these matchings.
     Moreover, it has no anticlique of order $4$.
     The complementary graph is a maximal $3$-fold $4$-colorable graph without a subgraph $K_4$.} 
    \end{center}
\end{figure}
}
$G$ is triangle free and has exactly $3$ perfect matchings, and deleting any single edge from $G$ produces a
triangle free graph with less than $3$ perfect matchings. Therefore, $\overline{G}$ has chromatic number $4$ and is maximal $3$-fold $4$-colorable.
However, $G$ has no anticlique of order $4$, so that $\omega(\overline{G}) < 4 = \chi(\overline{G})$.

The necessity part of the proof of Theorem \ref{TMaxambkcol} implies that a maximal ambiguously $k$-colorable graph $G$ with $\chi(G)=k$ has {\em exactly} two $k$-colorings.
This does not generalize to the higher dimensional case: Take two disjoint ambiguously $k$-colorable graphs $G_1,G_2$, with $\chi(G_1)=\chi(G_2)=k$, and
let $G:=G_1 * G_2$ be the graph obtained from their union by an edge $x_1 x_2$ for all $x_1 \in V(G_1)$ and all $x_2 \in V(G_2)$.
One readily checks that $\chi(G)=\chi(G_1)+\chi(G_2)=2k$, and that the $2k$-colorings of $G$ are obtained as a union of a $k$-coloring of $G_1$ and
a $k$-coloring of $G_2$. Hence $G$ has exactly four $2k$-colorings; by adding a single edge $e$ not yet present in $G$, the number of $2k$-colorings drops down to either two or zero,
as $e$ connects two nonadjacent vertices in $G_i$ for some $i \in \{1,2\}$ and either $G_i+e$ is uniquely $k$-colorable or not $k$-colorable at all.
Hence $G$ is maximal $4$-fold $2k$-colorable and, at the same time, maximal $3$-fold $2k$-colorable, with exactly $4$ $k$-colorings.
The construction generalizes in various obvious ways.

Concerning {\sc Hadwiger}'s Conjecture, I am optimistic: 

\begin{conjecture}
  Every maximal $d$-fold $k$-colorable graph is $(k-1)$-colorable or contains a clique minor of order $k$.
\end{conjecture}

{\bf Author's Address.}

{\sc Matthias Kriesell} \\
Ilmenau University of Technology \\
Weimarer Stra{\ss}e 25 \\
D--98693 Ilmenau \\
Germany

\end{document}